\documentclass[11pt,reqno]{amsart}
\usepackage{amsmath}
\usepackage{amssymb}
\usepackage[bookmarks]{hyperref}
\newtheorem{theorem}{Theorem}[section]
\newtheorem{lemma}[theorem]{Lemma}
\newtheorem{corollary}[theorem]{Corollary}
\newtheorem{proposition}[theorem]{Proposition}

\newtheorem{remark}[theorem]{Remark}
\newtheorem{definition}[theorem]{Definition}

\numberwithin{equation}{section}

\allowdisplaybreaks[1]
\begin{document}

\title[Commutants of Toeplitz operators]{The commutants of certain Toeplitz operators on weighted Bergman spaces}
\author{Trieu Le}
\address{Trieu Le, Department of Mathematics, University of Toronto, Toronto, Ontario, Canada M5S 2E4}
\email{trieu.le@utoronto.ca}
\begin{abstract}
For $\alpha>-1$, let $A^2_{\alpha}$ be the corresponding weighted Bergman space of the unit ball in $\mathbb{C}^n$. For a bounded measurable function $f$, let $T_f$ be the Toeplitz operator with symbol $f$ on $A^2_{\alpha}$. This paper describes all the functions $f$ for which $T_f$ commutes with a given $T_g$, where $g(z)=z_{1}^{L_1}\cdots z_{n}^{L_n}$ for strictly positive integers $L_1,\ldots, L_n$, or $g(z)=|z_1|^{s_1}\cdots |z_n|^{s_n}h(|z|)$ for non-negative real numbers $s_1,\ldots, s_n$ and a bounded measurable function $h$ on $[0,1)$.
\end{abstract}
\thanks{}
\date{}
\maketitle

\section{Introduction}
As usual, for any $z=(z_1,\ldots,z_n)\in\mathbb{C}^n$ we denote its Euclidean norm by $|z|$, which is $\sqrt{|z_1|^2+\cdots+|z_n|^2}$. Let $\mathbb{B}_n$ denote the open unit ball consisting of all $z\in\mathbb{C}^n$ with $|z|<1$. Let $\nu$ denote the Lebesgue measure on $\mathbb{B}_n$ normalized so that $\nu(\mathbb{B}_n)=1$. Fix a real number $\alpha>-1$. The weighted Lebesgue measure $\nu_{\alpha}$ on $\mathbb{B}_n$ is defined by $\mathrm{d}\nu_{\alpha}(z)=c_{\alpha}(1-|z|^2)^{\alpha}\mathrm{d}\nu(z)$, where $c_{\alpha}$ is a normalizing constant so that $\nu_{\alpha}(\mathbb{B}_n)=1$. A direct computation shows that $c_{\alpha}=\dfrac{\Gamma(n+\alpha+1)}{\Gamma(n+1)\Gamma(\alpha+1)}$. For $1\leq p\leq\infty$, let $L^{p}_{\alpha}$ denote the space $L^{p}(\mathbb{B}_n,\mathrm{d}\nu_{\alpha})$. Note that $L^{\infty}_{\alpha}$ is the same as $L^{\infty}=L^{\infty}(\mathbb{B}_n,\mathrm{d}\nu)$.

The weighted Bergman space $A^{2}_{\alpha}$ consists of all functions in $L^{2}_{\alpha}$ which are analytic on $\mathbb{B}_n$. It is well-known that $A^{2}_{\alpha}$ is a closed subspace of $L^{2}_{\alpha}$. We denote the inner product in $L^2_{\alpha}$ by $\langle \cdot, \cdot\rangle_{\alpha}$ and the corresponding norm by $\|\cdot\|_{2,\alpha}$.

For any multi-index $m=(m_1,\ldots, m_n)\in\mathbb{N}^n$ (here $\mathbb{N}$ denotes the set of all \emph{non-negative} integers), we write $|m|=m_1+\cdots+m_n$ and $m!=m_1!\cdots m_n!$. For any $z=(z_1,\ldots, z_n)\in\mathbb{C}^{n}$, we write $z^{m}=z_1^{m_1}\cdots z_n^{m_n}$ and $\bar{z}^{m}=\bar{z}_1^{m_1}\cdots \bar{z}_n^{m_n}$. The standard orthonormal basis for $A^{2}_{\alpha}$ is $\{e_{m}: m\in\mathbb{N}^n\}$, where
\begin{equation*}
e_{m}(z) = \Big[\dfrac{\Gamma(n+|m|+\alpha+1)}{m!\ \Gamma(n+\alpha+1)}\Big]^{1/2}z^{m},\ m\in\mathbb{N}^n, z\in\mathbb{B}_n.
\end{equation*}
For a more detailed discussion of $A^{2}_{\alpha}$, see Chapter 2 in \cite{Zhu2005}.

Since $A^2_{\alpha}$ is a closed subspace of the Hilbert space $L^2_{\alpha}$, there is an orthogonal projection $P_{\alpha}$ from $L^{2}_{\alpha}$ onto $A^{2}_{\alpha}$. For any function $f\in L^2_{\alpha}$ the Toeplitz operator with symbol $f$ is denoted by $T_f$, which is densely defined on $A^2_{\alpha}$ by $T_f\varphi = P_{\alpha}(f\varphi)$ for bounded analytic functions $\varphi$ on $\mathbb{B}_n$. If $f$ is a bounded function then $T_{f}$ is a bounded operator on $A^{2}_{\alpha}$ with $\|T_f\|\leq\|f\|_{\infty}$ and $(T_{f})^{*}=T_{\bar{f}}$. However there are unbounded functions $f$ that give rise to bounded operators $T_f$. If $f$ is an analytic function then $T_f$ is the multiplication operator on $A^2_{\alpha}$ with symbol $f$. The Toeplitz operator $T_f$ in this case is called \emph{analytic}. It is clear that if both $f$ and $g$ are bounded and analytic or conjugate analytic (that is, $\bar{f}$ and $\bar{g}$ are analytic) then $T_fT_g=T_gT_f$. Also if there are constants $a,b$ not both zero such that $af+bg$ is a constant function then it is clear that $T_f$ and $T_g$ commute. In the context of Toeplitz operators on the Hardy space of the unit circle, A. Brown and P. Halmos \cite{Brown1963} showed that these are the only cases where the operators $T_f$ and $T_g$ commute. For Toeplitz operators on the Bergman space of the unit disk, the situation becomes more complicated. The above Brown-Halmos result failed. In fact, if $f,g$ are radial functions, that is, $f(z)=f(|z|)$ and $g(z)=g(|z|)$ for almost all $z$, then both $T_f$ and $T_g$ are diagonal operators with respect to the standard orthonormal basis, hence, they commute. The problem that we are interested in is: if $T_f$ and $T_g$ commute on $A^2_{\alpha}$, what is the relation between the functions $f$ and $g$? Despite the difficulty of the general problem, several results have been known for Toeplitz operators on the Bergman space of the \emph{unit disk}:
\begin{enumerate}
\item If $g=z^N$ for some $N\geq 1$ then $f$ is analytic ({\v{Z}}. {\v{C}}u{\v{c}}kovi{\'c} \cite{Cuckovic1994}). This result was later extended to the case where $g$ is an arbitrary non-constant bounded analytic function by S. Axler, {\v{Z}}. {\v{C}}u{\v{c}}kovi{\'c} and N. Rao in \cite{Axler2000}.
\item If $f$ and $g$ are bounded harmonic functions, then either both functions are analytic or both are conjugate analytic or $af+bg$ is a constant for some constants $a$ and $b$ not both zero (Axler and {\v{C}}u{\v{c}}kovi{\'c} \cite{Axler1991}).
\item If $g$ is radial then $f$ is also radial ({\v{C}}u{\v{c}}kovi{\'c} and Rao \cite{Cuckovic1998}). In the same paper, they also characterized all bounded functions $f$ such that $T_f$ commutes with $T_g$ where $g(z)=z^{m_1}\bar{z}^{m_2}$ for integers $m_1,m_2\geq 0$.
\end{enumerate}

In this paper we generalize the results in (1) and (3) to Toeplitz operators on weighted Bergman spaces of the unit ball in higher dimensions. Let $\mathfrak{T}$ denote the $C^{*}-$algebra generated by $\{T_g: g\in L^{\infty}\}$. In addition to the result in (1), {\v{C}}u{\v{c}}kovi{\'c}  \cite{Cuckovic1994} showed that if $S$ is a bounded operator on the Bergman space of the unit disk such that $S$ belongs to $\mathfrak{T}$ and it commutes with $T_{z^N}$ for some integer $N\geq 1$, then $S=T_f$ for some bounded analytic function $f$. The following theorem generalizes this result.

\begin{theorem}\label{theorem-1} If $S$ is an operator in $\mathfrak{T}$ that commutes with $T_{z^{L_1}_{1}\cdots z^{L_n}_{n}}$ for some integers $L_1,\ldots, L_n\geq 1$, then there is a bounded analytic function $f$ on $\mathbb{B}_n$ so that $S=T_{f}$.
\end{theorem}

Due to the complicated setting of several variables, {\v{C}}u{\v{c}}kovi{\'c} and Rao's result in (3) no longer holds when $n\geq 2$. We will see later that if $f(z)=z_1\bar{z}_2$ then $T_f$ commutes with all $T_g$ whenever $g$ is a bounded radial function (a function $h$ on $\mathbb{B}_n$ is called a radial function if there is a function $\tilde{h}:[0,1)\rightarrow\mathbb{C}$ such that $h(z)=\tilde{h}(|z|)$ for almost all $z\in\mathbb{B}_n$). In fact, we will show that if $g$ is a non-constant bounded radial function then $T_f$ commutes with $T_g$ if and only if $f$ satisfies $f(\mathrm{e}^{\mathrm{i}\theta}z)=f(z)$ for almost all $z\in\mathbb{B}_n$ and almost all $\theta\in\mathbb{R}$.

Throughout the paper, an operator on $A^2_{\alpha}$ is said to be diagonal if it is diagonal with respect to the standard orthonormal basis of $A^2_{\alpha}$. In one dimension, all diagonal Toeplitz operators arise from radial functions. In higher dimensions, in order to get all diagonal operators we have to replace radial functions by functions that are invariant under the action of the $n$-torus on $\mathbb{B}_n$. More precisely, it was showed in \cite{Le-5} that for $f\in L^{\infty}$, the operator $T_f$ is diagonal if and only if $f(z_1,\ldots,z_n)=f(|z_1|,\ldots,|z_n|)$ for almost all $z\in\mathbb{B}_n$. Even though we are unable to describe all the functions $f\in L^2_{\alpha}$ such that $T_f$ commutes with a given non-trivial diagonal Toeplitz operator $T_g$, we have been successful in doing so when the function $g$ is of the from $g(z)=|z|^{2s_1}\cdots |z_n|^{2s_n}h(|z|)$, where $s_1,\ldots,s_n\geq 0$ and $h$ is a bounded function on $[0,1)$. The technique we use involves results about the zero sets of bounded analytic functions on the open unit disk. See Section \ref{section-3} for more detail.

Let $\mathcal{P}$ denote the space of all analytic polynomials in the variable $z=(z_1,\ldots,z_n)$. Then $\mathcal{P}$ is dense in $A^2_{\alpha}$. The following theorem is our second result in the paper.

\begin{theorem}\label{theorem-2} Let $g$ be a non-constant function in $\mathbb{B}_n$ such that for almost all $z\in\mathbb{B}_n$, $g(z)=|z_1|^{2s_1}\cdots |z_n|^{2s_n}h(|z|)$ with $h$ a bounded measurable function on $[0,1)$ and $s_1,\ldots, s_n\geq 0$. Then for $f\in L^2_{\alpha}$, $T_fT_g=T_gT_f$ on $\mathcal{P}$ if and only if $f(\mathrm{e}^{\mathrm{i}\theta}z)=f(z)$ for almost all $\theta\in\mathbb{R}$, almost all $z\in\mathbb{B}_n$, and for $1\leq j\leq n$ with $s_j\neq 0$, $f(z_1,\ldots,z_{j-1},|z_j|,z_{j+1},\ldots,z_n)=f(z)$ for almost all $z\in\mathbb{B}_n$.
\end{theorem}

\section{Commuting with analytic Toeplitz operators}\label{section-2}

The following result is well-known and its proof for the one dimensional case is in Proposition 7.2 in \cite{Zhu2007}. The proof for higher dimensional cases is similar. For the reader's convenience, we provide here the proof.

\begin{lemma}\label{lemma-2} Suppose $f=f_1+\bar{f}_2$, where $f_1,f_2\in A^2_{\alpha}$, such that $\|fp\|_{2,\alpha}\leq M\|p\|_{2,\alpha}$ for all analytic polynomials $p$, where $M$ is a fixed positive constant. Then $\|f\|_{\infty}\leq M$.
\end{lemma}
\begin{proof} For any $z\in\mathbb{B}_n$, let $k_{z}^{\alpha}(w)=(1-|z|^2)^{(n+\alpha+1)/2}(1-\langle w,z\rangle)^{-(n+\alpha+1)}$ for $w\in\mathbb{B}_n$. Then $k_{z}^{\alpha}$ is a normalized reproducing kernel for $A^2_{\alpha}$. This means $\|k_{z}^{\alpha}\|_{2,\alpha}=1$ and $\langle g,k_{z}^{\alpha}\rangle_{\alpha}=(1-|z|^2)^{(n+\alpha+1)/2}g(z)$ for all $g\in A^2_{\alpha}$. See Chapter 2 of \cite{Zhu2005} for more detail. Since $k_{z}^{\alpha}$ is analytic in a neighborhood of the closed unit ball, there is a sequence $\{p_{s}\}_{s=1}^{\infty}$ of analytic polynomials converging uniformly to $k^{\alpha}_{z}$ on $\bar{\mathbb{B}}_n$. It then follows that $\lim\limits_{s\rightarrow\infty}\|f p_s - f k_{z}^{\alpha}\|_{2,\alpha}\rightarrow 0$. Hence
$
\langle fk_{z}^{\alpha}, k_{z}^{\alpha}\rangle_{\alpha} = \lim\limits_{s\rightarrow\infty}\langle fp_s, k_{z}^{\alpha}\rangle_{\alpha}.
$ Now for any integer $s\geq 1$, we have $|\langle fp_s, k_{z}^{\alpha}\rangle_{\alpha}|\leq \|fp_s\|_{2,\alpha}\leq M\|p_s\|_{2,\alpha}$. So we conclude that $
|\langle fk_{z}^{\alpha}, k_{z}^{\alpha}\rangle_{\alpha}|\leq \lim\limits_{s\rightarrow\infty} M\|p_s\|_{2,\alpha} = M\|k_{z}^{\alpha}\|_{2,\alpha}=M$. On the other hand,
\begin{align*}
\langle fk_{z}^{\alpha}, k_{z}^{\alpha}\rangle_{\alpha} & = \langle f_1 k_{z}^{\alpha}, k_{z}^{\alpha}\rangle_{\alpha} + \langle\bar{f}_2 k_{z}^{\alpha}, k_{z}^{\alpha}\rangle_{\alpha}\\
& = (1-|z|^2)^{(n+1+\alpha)/2}f_1(z)k_{z}^{\alpha}(z) + (1-|z|^2)^{(n+1+\alpha)/2}\bar{f}_2(z)k_{z}^{\alpha}(z)\\
& = f_1(z)+\bar{f}_2(z) = f(z).
\end{align*}
So from the above inequality, $|f(z)|\leq M$ for any $z\in\mathbb{B}_n$. This shows that $\|f\|_{\infty}\leq M$.
\end{proof}

For any multi-indexes $m,k\in\mathbb{N}^n$, there is a positive real number $d(m,k)$ such that $e_m e_k = d(m,k)e_{m+k}$. Strictly speaking, $d(m,k)$ must be written as $d_{n,\alpha}(m,k)$ because it depends also on $n$ and $\alpha$. But to simplify the notation and since $n$ and $\alpha$ are fixed throughout the paper, we drop the sub-indexes.

It is immediate that $d(m,k)=d(k,m)$ and $d(0,k)=d(m,0)=1$ for $m,k\in\mathbb{N}^n$. Also for any $m,k,l\in\mathbb{N}^n$, since $(e_{m}e_{k})e_{l}=e_{m}(e_{k}e_{l})$, we have
\begin{equation}\label{eqn-403}
d(m,k)d(m+k,l) = d(m,k+l)d(k,l).
\end{equation}
Using the explicit formulas for $e_{m}$ and $e_{k}$, we obtain
\begin{align}\label{eqn-408}
(d(m,k))^2 & = \dfrac{\Gamma(n+|m|+\alpha+1)\ \Gamma(n+|k|+\alpha+1)}{\Gamma(n+|m|+|k|+\alpha+1)\ \Gamma(n+\alpha+1)}\dfrac{(m+k)!}{m!\ k!}.
\end{align}

The following lemma characterizes analytic Toeplitz operators on $A^2_{\alpha}$ in terms of its matrix with respect to the standard orthonormal basis. Even though the matrix of an analytic Toeplitz operator is not the usual \emph{analytic Toeplitz matrix}, it becomes one after scaling each matrix entry by a factor depending on the position of the entry.

For an $n-$tuple of integers $r=(r_1,\ldots,r_n)\in\mathbb{Z}^n$ we write $r\succeq 0$ if $r_1,\ldots,r_n\geq 0$ and write $r\nsucceq 0$ if otherwise. For $m,k\in\mathbb{N}^n$ we write $m\succeq k$ (respectively, $m\nsucceq k$) if $m-k\succeq 0$ (respectively, $m-k\nsucceq 0$).

\begin{lemma}\label{lemma-3} Suppose $S$ is a linear operator (not necessarily bounded) on $A^2_{\alpha}$ whose domain contains the space $\mathcal{P}$ of all analytic polynomials. Then there is a function $f\in A^2_{\alpha}$ such that $T_f = S$ on $\mathcal{P}$ if and only if $\langle S e_{m}, e_{k}\rangle_{\alpha}=0$ whenever $k\nsucceq m$ and for any $l\in\mathbb{N}^n$, $\dfrac{1}{d(l,m)}\langle S e_m, e_{m+l}\rangle_{\alpha}$ is independent of $m\in\mathbb{N}^n$.
\end{lemma}

\begin{proof}
Suppose $f\in A^2_{\alpha}$ and it has the expansion $f=\sum\limits_{l\in\mathbb{N}^n}a_{l}e_{l}$. Then for any $m$ in $\mathbb{N}^n$, we have
\begin{equation*}
T_f e_{m} = fe_{m} = \sum\limits_{l\in\mathbb{N}^n}a_{l}e_{l}e_{m} = \sum\limits_{l\in\mathbb{N}^n} a_{l}d(l,m)e_{m+l}.
\end{equation*}
Therefore, $
\langle T_f e_{m}, e_{k}\rangle_{\alpha}  = \begin{cases}
0 & \text{ if } k\nsucceq m\\
a_{l} d(l,m) & \text{ if } k=m+l.
\end{cases}$

This shows that if $S = T_f$ on $\mathcal{P}$ then $\langle S e_{m}, e_{k}\rangle_{\alpha}=0$ whenever $k\nsucceq m$ and for any $l\in\mathbb{N}^n$, $\dfrac{1}{d(l,m)}\langle Te_m, e_{m+l}\rangle_{\alpha}=a_{l}$, which is independent of $m$.

Now suppose $S$ has the above property. Let $a_{l}=\langle Te_{0}, e_{l}\rangle_{\alpha}$ for each $l$ in $\mathbb{N}^n$. Then by assumption, $a_{l}=\dfrac{1}{d(l,m)}\langle S e_m, e_{m+l}\rangle_{\alpha}$ for all $m$ in $\mathbb{N}^n$. We have
\begin{equation*}
\sum\limits_{l\in\mathbb{N}^n}|a_{l}|^2 = \sum\limits_{l\in\mathbb{N}^n}|\langle S e_{0}, e_{l}\rangle_{\alpha}|^2 = \|S e_{0}\|^2_{2,\alpha}<\infty.
\end{equation*}
So the function $f=\sum\limits_{l\in\mathbb{N}^n}a_{l}e_{l}$ is an element of $A^2_{\alpha}$ and we have $\langle T_f e_{m}, e_{k}\rangle_{\alpha}=\langle S e_{m}, e_{k}\rangle_{\alpha}$ for all $m,k$ in $\mathbb{N}^n$. Hence $T_f=S$ on $\mathcal{P}$.
\end{proof}

Suppose $l\in\mathbb{N}^n$ is a multi-index. For any $m,k\in\mathbb{N}^n$ we have
\begin{align*}
\langle T_{e_{l}}e_{m}, e_{k}\rangle_{\alpha} & = d(l,m)\langle e_{l+m}, e_{k}\rangle_{\alpha}
= \begin{cases}
0 & \text{ if } k\neq l+m,\\
d(l,m) & \text{ if } k=l+m.
\end{cases}
\end{align*}
This implies $T_{e_{l}}e_{m}=d(l,m)e_{l+m}$ and
$
T_{\bar{e}_{l}}e_{k} = \begin{cases}
0 & \text{if } k\nsucceq l,\\
d(l,k-l)e_{k-l} & \text{if } k\succeq l.
\end{cases}
$

For any linear operator $S$ on $A^2_{\alpha}$ whose domain contains the space of all analytic polynomials, we have
\begin{align*}
\langle [S,T_{e_{l}}]e_{m},e_{k}\rangle_{\alpha} & = \langle ST_{e_{l}}e_{m},e_{k}\rangle_{\alpha}-\langle Se_{m},T_{\bar{e}_{l}}e_{k}\rangle_{\alpha}\notag\\
& = d(l,m)\langle S e_{m+l},e_{k}\rangle_{\alpha} - \begin{cases}
0 & \text{if } k\nsucceq l,\\
d(l,k-l)e_{k-l} & \text{if } k\succeq l.
\end{cases}
\end{align*}
This shows that for $m,k\in\mathbb{N}^n$,
\begin{align}
\langle [S,T_{e_{l}}]e_{m},e_{k}\rangle_{\alpha} & = d(l,m)\langle S e_{m+l}, e_{k}\rangle_{\alpha} \text{ if } k\nsucceq l,\label{eqn-409}\\
\text{and }\langle [S,T_{e_{l}}]e_{m},e_{k+l}\rangle_{\alpha} & = d(l,m)\langle S e_{m+l}, e_{k+l}\rangle_{\alpha}-d(l,k)\langle S e_{m}, e_{k}\rangle_{\alpha}.\label{eqn-410}
\end{align}

\begin{lemma}\label{lemma-4} Suppose that $S$ is an operator (not necessarily bounded) on $A^2_{\alpha}$ whose domain contains $\mathcal{P}$ and that $S$ commutes with $T_{e_{L}}$ where $L=(L_1,\ldots,L_n)\in\mathbb{N}^n$ with $L_1,\ldots,L_n\geq 1$. Suppose $l\in\mathbb{N}^n$ such that the operator $K=[S,T_{e_{l}}]$ is a compact operator on $A^2_{\alpha}$. Then for any $m,k\in\mathbb{N}^n$ we have $
\dfrac{d(l,m)}{d(l,k)}\langle S e_{m+l}, e_{k+l}\rangle_{\alpha} = \langle S e_{m}, e_{k}\rangle_{\alpha}.$
\end{lemma}

\begin{proof}
Since $S$ commutes with $T_{e_{L}}$, it commutes with $T_{e_{sL}}$ for any positive integer $s$ because $T_{e_{sL}}$ is a multiple of $(T_{e_{L}})^{s}$. Then \eqref{eqn-410} implies that for any $m,k\in\mathbb{N}^n$ and $s\in\mathbb{N}$,
\begin{equation*}
\dfrac{d(sL,m)}{d(sL,k)}\langle S e_{m+sL}, e_{k+sL}\rangle_{\alpha} = \langle S e_{m}, e_{k}\rangle_{\alpha}.
\end{equation*}
Now for any $m,k\in\mathbb{N}^n$, and $s\in\mathbb{N}$,
\begin{align*}
& \dfrac{d(l,m)}{d(l,k)}\langle S e_{m+l},e_{k+l}\rangle_{\alpha} - \langle S e_{m},e_{k}\rangle_{\alpha}\\
& = \dfrac{d(l,m)}{d(l,k)}\dfrac{d(sL,m+l)}{d(sL,k+l)}\langle S e_{m+l+sL}, e_{k+l+sL}\rangle_{\alpha}-\dfrac{d(sL,m)}{d(sL,k)}\langle S e_{m+sL},e_{k+sL}\rangle_{\alpha}\\
& = \dfrac{d(sL,m)d(sL+m,l)}{d(sL,k)d(sL+k,l)}\langle S e_{m+l+sL}, e_{k+l+sL}\rangle_{\alpha} -\dfrac{d(sL,m)}{d(sL,k)}\langle S e_{m+sL},e_{k+sL}\rangle_{\alpha}\\
& = \dfrac{d(sL,m)}{d(sL,k)d(sL+k,l)}\Big\{d(sL+m,l)\langle S e_{m+l+sL}, e_{k+l+sL}\rangle_{\alpha}\\
&\phantom{\dfrac{d(sL,m)}{d(sL,k)d(sL+k,l)}}\quad\quad - d(sL+k,l)\langle S e_{m+sL},e_{k+sL}\rangle_{\alpha}\Big\}\\
& = \dfrac{d(sL,m)}{d(sL,k)d(sL+k,l)}\langle [S,T_{e_{l}}]e_{sL+m}, e_{sL+k+l}\rangle_{\alpha}\quad\quad\text{(by \eqref{eqn-410})}\\
& = \dfrac{d(sL,m)}{d(sL,k)d(sL+k,l)}\langle K e_{sL+m}, e_{sL+k+l}\rangle_{\alpha}.
\end{align*}
Now using \eqref{eqn-408} we have
\begin{align*}
&\Big[\dfrac{d(sL,m)}{d(sL,k)d(sL+k,l)}\Big]^2\\
& = \dfrac{\Gamma(n+s|L|+\alpha+1)\ \Gamma(n+|m|+\alpha+1)}{\Gamma(n+s|L|+|m|+\alpha+1)\ \Gamma(n+\alpha+1)}\dfrac{(sL+m)!}{(sL)!\ m!}\\
&\quad\times \dfrac{\Gamma(n+s|L|+|k|+\alpha+1)\ \Gamma(n+\alpha+1)}{\Gamma(n+s|L|+\alpha+1)\ \Gamma(n+|k|+\alpha+1)}\dfrac{(sL)!\ k!}{(sL+k)!}\\
&\quad\times \dfrac{\Gamma(n+s|L|+|k|+|l|+\alpha+1)\ \Gamma(n+\alpha+1)}{\Gamma(n+s|L|+|k|+\alpha+1)\ \Gamma(n+|l|+\alpha+1)}\dfrac{(sL+k)!\ l!}{(sL+k+l)!}\\
& = C(n,\alpha,m,k,l,L)\dfrac{\Gamma(n+s|L|+|l|+|k|+\alpha+1)}{\Gamma(n+s|L|+|m|+\alpha+1)}\dfrac{(sL+m)!}{(sL+k+l)!}\\
& \approx C(n,\alpha,m,k,l,L)(s|L|)^{|l|+|k|-|m|}\prod_{j=1}^{n}(sL_{j})^{m_j-k_j-l_j}\\
&\quad\text{(by Stirling's formula for the Gamma function)}\\
& = \tilde{C}(n,\alpha,m,k,l,L).
\end{align*}
This shows that $\dfrac{d(sL,m)}{d(sL,k)d(sL+k,l)}$ is bounded when $s\rightarrow\infty$. On the other hand, $\lim\limits_{s\rightarrow\infty}\langle K e_{sL+m}, e_{sL+k+l}\rangle_{\alpha}=0$ because $K$ is compact. So we conclude that $
\dfrac{d(l,m)}{d(l,k)}\langle S e_{m+l},e_{k+l}\rangle_{\alpha} = \langle S e_{m},e_{k}\rangle_{\alpha}$ for all $m,k\in\mathbb{N}^n$.
\end{proof}

\begin{proof}[Proof of Theorem \ref{theorem-1}] Put $L=(L_1,\ldots,L_n)$. Then by assumption, $S$ commutes with $T_{e_{L}}$. Since $S$ belongs to $\mathfrak{T}$, it is well-known that $[S,T_{e_{l}}]$ is compact for all $l\in\mathbb{N}$ (see \cite{Coburn1973} for more detail). Now Lemma \ref{lemma-4} shows that
\begin{equation*}
\dfrac{d(l,m)}{d(l,k)}\langle S e_{m+l}, e_{k+l}\rangle_{\alpha} = \langle S e_{m}, e_{k}\rangle_{\alpha}\quad\text{ for any } m,k,l\in\mathbb{N}^n.
\end{equation*}
Put $m=0$ we see that for each $k\in\mathbb{N}^n$, $\dfrac{1}{d(l,k)}\langle S e_{l}, e_{k+l}\rangle_{\alpha} = \langle S e_{0}, e_{k}\rangle_{\alpha}$, which is independent of $l\in\mathbb{N}^n$.

Now suppose $m,k\in\mathbb{N}^n$ such that $k\nsucceq m$. Then there is an integer $1\leq j\leq n$ so that $k_j<m_j$. Consider first the case $k_j=0$. Put $l=(L_1,\ldots, L_{j-1}, L_j-1, L_{j+1},\ldots, L_{n})$. Then we have $k+l\nsucceq L$ but $m+l\geq L$. Hence,
\begin{align*}
\langle S e_{m}, e_{k}\rangle_{\alpha} & = \dfrac{d(l,m)}{d(l,k)}\langle S e_{m+l}, e_{k+l}\rangle_{\alpha}\\
& =\dfrac{d(l,m)}{d(l,k)d(L,m+l-L)}\langle S T_{e_{L}}e_{m+l-L}, e_{k+l}\rangle_{\alpha}\\
& =\dfrac{d(l,m)}{d(l,k)d(L,m+l-L)}\langle T_{e_{L}}Se_{m+l-L}, e_{k+l}\rangle_{\alpha}\\
& = \dfrac{d(l,m)}{d(l,k)d(L,m+l-L)}\langle Se_{m+l-L}, T_{\bar{e}_{L}}e_{k+l}\rangle_{\alpha}=0
\end{align*}
since $T_{\bar{e}_{L}}e_{k+l}=0$.

Now consider the case $k_j>0$. Let $\tilde{k}=k-k_j\delta_j$ and $\tilde{m}=m-k_j\delta_j$, where $\delta_{j}=(\delta_{1j},\ldots,\delta_{nj})$. Then $\tilde{k},\tilde{m}$ are in $\mathbb{N}^n$ and $0=\tilde{k}_j<\tilde{m}_j$. We have
\begin{align*}
\langle S e_{m}, e_{k}\rangle_{\alpha} & = \langle S e_{\tilde{m}+m_j\delta_j}, e_{\tilde{k}+m_j\delta_j}\rangle_{\alpha}\\
& = \dfrac{d(m_j\delta_j,\tilde{k})}{d(m_j\delta_{j},\tilde{m})}\langle S e_{\tilde{m}}, e_{\tilde{k}}\rangle_{\alpha}\\
& = 0\text{ (by the case considered above).}
\end{align*}
By Lemma \ref{lemma-3} there is a function $f\in A^2_{\alpha}$ so that $S=T_f$ on the space of analytic polynomials. Since $S$ is a bounded operator, Lemma \ref{lemma-2} implies that $f$ is bounded and $\|f\|_{\infty}\leq\|S\|$. Consequently, $S=T_f$ on $A^2_{\alpha}$ because they are bounded operators that agree on a dense subset of $A^2_{\alpha}$.
\end{proof}

We now discuss the necessity of the condition $L_1,\ldots, L_n\geq 1$ in Theorem \ref{theorem-1}. When $n=1$ this condition is necessary because $\mathfrak{T}$ obviously contains non-analytic Toeplitz operators that commute with $T_{e_{0}}\equiv I$. For the case $n\geq 2$, we will show that there is an operator $S$ in $\mathfrak{T}$ such that $S$ commutes with $T_{z_1},\ldots, T_{z_{n-1}}$ but it does not commute with $T_{z_n}$. Recall that for $1\leq j\leq n$, $\delta_{j}$ denotes $(\delta_{1j},\ldots,\delta_{nj})$.

\begin{proposition}\label{prop-4} Suppose $n\geq 2$. For any $\varphi\in A^{2}_{\alpha}$, define
\begin{equation}\label{eqn-411}
S\varphi = S\Big(\sum\limits_{m\in\mathbb{N}^n}\langle\varphi, e_{m}\rangle_{\alpha}e_{m}\Big) = \sum_{\substack{m\in\mathbb{N}^n\\
m_{n}=0
}}d(m,\delta_n)\langle\varphi, e_{m}\rangle_{\alpha} e_{m+\delta_n}.
\end{equation}
Then the following statements hold true:
\begin{enumerate}
\item $S$ is a compact operator on $A^2_{\alpha}$ and hence it belongs to $\mathfrak{T}$,
\item $S$ commutes with $T_{z_{1}},\ldots, T_{z_{n-1}}$,
\item $S$ does not commute with $T_{z_{n}}$ and hence $S$ is not an analytic Toeplitz operator.
\end{enumerate}
\end{proposition}

\begin{proof}
From the definition \eqref{eqn-411} of $S$, we see that $Se_{m}=d(m,\delta_n)e_{m+\delta_n}$ if $m_n=0$ and $Se_{m}=0$ if $m_n>0$. For any $m\in\mathbb{N}^n$ with $m_n=0$, formula \eqref{eqn-408} gives
\begin{align*}
(d(m,\delta_n))^2 & = \dfrac{\Gamma(n+|m|+\alpha+1)\ \Gamma(n+|\delta_n|+\alpha+1)}{\Gamma(n+|m|+|\delta_n|+\alpha+1)\ \Gamma(n+\alpha+1)}\dfrac{(m+\delta_n)!}{m!\delta_n!}\\
& = \dfrac{\Gamma(n+|m|+\alpha+1)\ \Gamma(n+\alpha+2)}{\Gamma(n+|m|+\alpha+2)\ \Gamma(n+\alpha+1)}\dfrac{(m+\delta_n)!}{m!}\\
& = \dfrac{(n+\alpha+1)(m_n+1)}{n+|m|+\alpha+1}\\
& = \dfrac{n+\alpha+1}{n+|m|+\alpha+1}\text{ (since $m_n=0$)}.
\end{align*}
Thus, $\lim\limits_{\substack{|m|\rightarrow\infty\\ m_n=0}}d(m,\delta_n)=0$. This shows that the operator $S$ is not only bounded but also compact on $A^2_{\alpha}$. On the other hand, it is well-known that $\mathfrak{T}$ contains the ideal of compact operators on $A^2_{\alpha}$ (see \cite{Coburn1973}). Hence $S$ belongs to $\mathfrak{T}$.

Now let $j$ be an integer in $\{1,\ldots,n-1\}$. For $m\in\mathbb{N}^n$ we have
\begin{align*}
ST_{e_{\delta_j}} e_{m} & = S\big(d(m,\delta_j)e_{m+\delta_j}\big)\\
& = \begin{cases}
0 & \text{if } m_n>0\\
d(m,\delta_j)d(m+\delta_j,\delta_n)e_{m+\delta_j+\delta_n} & \text{if } m_n=0,
\end{cases}\\
T_{e_{\delta_j}}S e_{m} & = \begin{cases}
0 & \text{if } m_n>0\\
T_{e_{\delta_j}}\big(d(m,\delta_n)e_{m+\delta_n}\big) & \text{if } m_n=0
\end{cases}\\
& = \begin{cases}
0 & \text{if } m_n>0\\
d(m,\delta_n)d(\delta_j,m+\delta_n)e_{m+\delta_n+\delta_j} & \text{if } m_n=0
\end{cases}\\
& = \begin{cases}
0 & \text{if } m_n>0\\
d(m,\delta_j)d(m+\delta_j,\delta_n)e_{m+\delta_n+\delta_j} & \text{if } m_n=0.
\end{cases}
\end{align*}
Thus $ST_{e_{\delta_j}}e_{m}=T_{e_{\delta_j}}Se_{m}$ for all $m\in\mathbb{N}^n$. This shows that $S$ commutes with $T_{e_{\delta_j}}$ (hence $T_{z_j}$) for $1\leq j\leq n-1$. Now for $m\in\mathbb{N}^n$ with $m_n=0$, we have
\begin{align*}
ST_{e_{\delta_n}} e_{m} & = S\big(d(\delta_n,m)e_{m+\delta_n}\big) = 0,\\
T_{e_{\delta_n}}S e_{m} & = T_{e_{\delta_{n}}}\big(d(m,\delta_n)e_{m+\delta_n}\big) = d(\delta_n,m+\delta_n)d(m,\delta_n)e_{m+2\delta_{n}}\neq 0.
\end{align*}
This shows that $ST_{e_{\delta_n}}\neq T_{e_{\delta_n}}S$, so $S$ does not commute with $T_{e_{\delta_n}}$. Since $T_{z_n}$ is a nonzero multiple of $T_{e_{\delta_n}}$, $S$ does not commute with $T_{z_n}$ either.
\end{proof}

\section{Commuting with diagonal Toeplitz operators}\label{section-3}

In the first part of this section we use results from complex analysis of one variable, more precisely, results about zeros of bounded analytic functions on the open unit disk, to obtain some function-theoretic results which are crucial for the proof of Theorem \ref{theorem-2}. The proof of Theorem \ref{theorem-2} itself will be presented at the end of the section.

For the rest of the paper, $N_{*}$ denotes the set of all positive integers. For any $1\leq j\leq n$, let $\sigma_{j}:\mathbb{N}_{*}\times\mathbb{N}_{*}^{n-1}\longrightarrow\mathbb{N}_{*}^{n}$ be the map defined by the formula $\sigma_{j}(s,(r_1,\ldots,r_{n-1}))=(r_1,\ldots,r_{j-1},s,r_{j},\ldots,r_{n-1})$ for all $s\in\mathbb{N}_{*}$ and $(r_1,\ldots,r_{n-1})\in\mathbb{N}_{*}^{n-1}$. If $M$ is a subset of $\mathbb{N}_{*}^n$ and $1\leq j\leq n$, we define
\begin{equation*}\label{eqn-2}
\widetilde{M}_{j} = \Big\{\tilde{r}=(r_1,\ldots,r_{n-1})\in\mathbb{N}_{*}^{n-1}: \sum\limits_{\substack{s\in\mathbb{N}_{*}\\
                                                           \sigma_{j}(s,\tilde{r})\in M}}
                                            \dfrac{1}{s+1}=\infty\Big\}.
\end{equation*}

\begin{definition}
We say that $M$ has property (P) if one of the following statements holds.
\begin{enumerate}
\item $M=\emptyset$, or
\item $M\neq\emptyset$, $n=1$ and $\sum\limits_{s\in M}\dfrac{1}{s+1}<\infty$, or
\item $M\neq\emptyset$, $n\geq 2$ and for any $1\leq j\leq n$, the set $\widetilde{M}_{j}$ has property (P) as a subset of $\mathbb{N}_{*}^{n-1}$.
\end{enumerate}
\end{definition}

The following observations are then immediate.
\begin{enumerate}
\item If $M\subset\mathbb{N}_{*}$ and $M$ does not have property (P) then $\sum_{s\in M}\frac{1}{s+1}=\infty$. If $M\subset\mathbb{N}_{*}^{n}$ with $n\geq 2$ and $M$ does not have property (P) then $\widetilde{M}_{j}$ does not have property (P) as a subset of $\mathbb{N}_{*}^{n-1}$ for some $1\leq j\leq n$.

\item If $M_1$ and $M_2$ are subsets of $\mathbb{N}_{*}^{n}$ that both have property (P) then $M_1\cup M_2$ also has property (P).

\item If $M\subset\mathbb{N}_{*}^n$ has property (P) and $l\in\mathbb{Z}^{n}$ then $(M+l)\cap\mathbb{N}_{*}^n$ also has property (P). Here $M+l=\{m+l: m\in M\}$.

\item If $M\subset\mathbb{R}^n$ has property (P) then $\mathbb{N}_{*}\times M$ also has property (P) as a subset of $\mathbb{N}_{*}^{n+1}$. This can be showed by induction on $n$.

\item The set $\mathbb{N}_{*}^n$ does not have property (P) for all $n\geq 1$. This together with (2) shows that if $M\subset\mathbb{N}_{*}^n$ has property (P) then $\mathbb{N}_{*}^n\backslash M$ does not have property (P).

\end{enumerate}

\begin{proposition}\label{prop-1} Let $\mathbb{K}$ denote the right half of the complex plane. Let $F:\mathbb{K}^{n}\rightarrow\mathbb{C}$ be an analytic function. Suppose there exists a polynomial $p$ such that $|F(z)|\leq p(|z|)$ for all $z\in\mathbb{K}^n$. Put $Z(F)=\{r\in\mathbb{N}_{*}^{n}: F(r)=0\}$. If $Z(F)$ does not have property (P), then $F$ is identically zero in $\mathbb{K}^n$.
\end{proposition}
\begin{proof}
Suppose $F$ is a function that satisfies the hypothesis of the proposition and that $Z(F)$ does not have property (P). We will show that $F(z)=0$ for all $z\in\mathbb{K}^n$ by induction on $n$.

Consider the case $n=1$. Write $p(|z|)=a_{0}+\cdots + a_{d}|z|^d$ for some positive integer $d$. For $z\in\mathbb{K}$, since $\max\{|z|,1\}\leq |z+1|$, we have $p(|z|)\leq (|a_0|+\cdots+|a_d|)|z+1|^d$. Let $G(z)=F(z)/(z+1)^{d}$ for $z\in\mathbb{K}$. Then $G$ is a bounded analytic function on $\mathbb{K}$ and $Z(G)=Z(F)$. Now define
\begin{equation*} H(z) = G\Big(\dfrac{1+z}{1-z}\Big)\quad\quad (|z|<1).\end{equation*}
Then $H$ is a bounded analytic function on the unit disk. We have $H(\theta)=0$ for all $\theta=\frac{r-1}{r+1}$ with $r\in Z(G)$. Since
\begin{equation*}
\sum\limits_{H(\theta)=0}(1-|\theta|)\geq\sum\limits_{r\in Z(G)}\Big(1-\big|\dfrac{r-1}{r+1}\big|\Big) = \sum\limits_{r\in Z(G)}\dfrac{2}{r+1} = \sum\limits_{r\in Z(F)}\dfrac{2}{r+1} =\infty,
\end{equation*}
Corollary to Theorem 15.23 in \cite{Rudin1987} shows that $H$ is identically zero on the unit disk. Thus $G$ is identically zero in $\mathbb{K}$, which implies that $F$ is identically zero in $\mathbb{K}$.

Now suppose that the conclusion of the proposition holds whenever $n\leq N$ for some integer $N\geq 1$. Consider the case $n=N+1$. Since $Z(F)$ does not have property (P), $\widetilde{Z(F)}_{j}$ does not have property (P) for some $1\leq j\leq N+1$. Without loss of generality, we may assume that $j=N+1$. For any $\tilde{r}$ in $\widetilde{Z(F)}_{N+1}$, put $M_{\tilde{r}}=\{s\in\mathbb{N}_{*}:(\tilde{r},s)\in Z(F)\}$. Then $\sum\limits_{s\in M_{\tilde{r}}}\frac{1}{s+1}=\infty$. Put $F_{\tilde{r}}(\zeta)=F(\tilde{r},\zeta)$ for $\zeta\in\mathbb{K}$. Then $F_{\tilde{r}}$ is analytic in $\mathbb{K}$ with $Z(F_{\tilde{r}})=M_{\tilde{r}}$ and $|F_{\tilde{r}}(\zeta)|\leq p(|(\tilde{r},\zeta)|)$ for all $\zeta\in\mathbb{K}$. Since $M_{\tilde{r}}$ does not have property (P), the proposition in the case $n=1$ implies that $F_{\tilde{r}}(\zeta)=0$ for all $\zeta\in\mathbb{K}$. Hence we have $F(\tilde{r},\zeta)=0$ for all $\zeta\in\mathbb{K}$ and all $\tilde{r}\in\widetilde{Z(F)}_{N+1}$. But $\widetilde{Z(F)}_{N+1}$ does not have property (P), the induction hypothesis shows that $F(\tilde{z},\zeta)=0$ for all $\zeta\in\mathbb{K}$ and all $\tilde{z}\in\mathbb{K}^{N}$. Thus $F$ is identically zero on $\mathbb{K}^{N+1}$.
\end{proof}

\begin{lemma}\label{lemma-1} For any function $f\in L^{1}(\mathbb{B}_n,\mathrm{d}\nu)$ and any $l\in\mathbb{Z}^n$, put
\begin{equation*}
Z(f,l) = \{m\in\mathbb{N}^n: m+l\succeq 0\text{ and }\int\limits_{\mathbb{B}_n}f(z)z^{m+l}\bar{z}^{m}\mathrm{d}\nu=0\}.
\end{equation*}
If $Z(f,l)$ does not have property (P) then it is the set of all $m\in\mathbb{N}^n$ with $m+l\succeq 0$.
\end{lemma}
\begin{proof}
Suppose $l=(l_1,\ldots,l_n)$ where $l_j\in\mathbb{Z}$ for $j=1,\ldots,n$. Put $l^{*}=(|l_1|,\ldots,|l_n|)$, $l^{+}=\frac{1}{2}(l^{*}+l)$ and $l^{-}=\frac{1}{2}(l^{*}-l)$. Then $l^{+}, l^{-}\succeq 0$ and $l=l^{+}-l^{-}$. Also, for any $m\in\mathbb{N}^n$, we have $m+l\succeq 0$ if and only if $m-l^{-}\succeq 0$. For $m\in Z(f,l)$, put $k=m-l^{-}$, then $k\succeq 0$ and
\begin{align}\label{eqn-1}
0  = \int\limits_{\mathbb{B}_n}f(z)z^{k+l^{+}}\bar{z}^{k+l^{-}}\mathrm{d}\nu & = \int\limits_{\mathbb{B}_n}f(z)z^{l^{+}}\bar{z}^{l^{-}}|z_1|^{2k_1}\cdots|z_n|^{2k_n}\mathrm{d}\nu\\
& = \int\limits_{\substack{\mathbb{B}_n\\ |z_1|,\ldots,|z_n|>0}}f(z)z^{l^{+}}\bar{z}^{l^{-}}|z_1|^{2k_1}\cdots|z_n|^{2k_n}\mathrm{d}\nu.
\end{align}

For any $\zeta=(\zeta_1,\ldots,\zeta_n)\in\bar{\mathbb{K}}^n$, define
\begin{equation*}
F(\zeta) = \int\limits_{\substack{\mathbb{B}_n\\ |z_1|,\ldots,|z_n|>0}}f(z)z^{l^{+}}\bar{z}^{l^{-}}|z_1|^{2\zeta_1}\cdots|z_n|^{2\zeta_n}\mathrm{d}\nu.
\end{equation*}
Here for a complex number $w$ and a real number $t>0$, $t^{w}=\exp(w\log t)$, where $\log$ is the principle branch of the logarithmic function. Since $|t^{w}|\leq 1$ for all $0<t<1$ and $w\in\mathbb{C}$ with $\Re(w)\geq 0$, the function $F$ is well-defined, bounded, and in fact continuous on $\bar{\mathbb{K}}^n$. Now an application of Morera's Theorem shows that $F$ is analytic on $\mathbb{K}^n$.

Next, \eqref{eqn-1} shows that $F(k)=0$ for all $k\in Z(f,l)-l^{-}$. Since $Z(f,l)$ does not have property (P), $Z(F)=\{r\in\mathbb{N}_{*}^n: F(r)=0\}$ does not have property (P) either. Proposition \ref{prop-1} and the continuity of $F$ on $\bar{\mathbb{K}}^n$ now imply that $F(\zeta)=0$ for all $\zeta\in\bar{\mathbb{K}}^n$. In particular, \eqref{eqn-1} holds for all $k\in\mathbb{N}^n$. The conclusion of the lemma then follows.
\end{proof}

\begin{corollary}\label{cor-1} Suppose $f\in L^{1}(\mathbb{B}_n,\mathrm{d}\nu)$ such that for all $l\in\mathbb{Z}^n$ the set $Z(f,l)$ (as in Lemma \ref{lemma-1}) does not have property (P). Then $f(z)=0$ for almost all $z\in\mathbb{B}_n$.
\end{corollary}
\begin{proof}
Lemma \ref{lemma-1} shows that $Z(f,l)=\mathbb{N}^n\cap(\mathbb{N}^n-l)$ for all $l\in\mathbb{Z}^n$. This implies that $\displaystyle\int\limits_{\mathbb{B}_n}f(z)z^{m}\bar{z}^k\mathrm{d}\nu=0$ for all $m,k\in\mathbb{N}^n$. Since the span of $\{z^m\bar{z}^k: m,k\in\mathbb{N}^n\}$ is dense in $C(\bar{\mathbb{B}}_n)$ we conclude that $f(z)=0$ for almost all $z\in\mathbb{B}_n$.
\end{proof}

\begin{corollary}\label{cor-2} Let $\gamma=(\gamma_1,\ldots,\gamma_n)$ be an $n$-tuple of integers and let $f$ be in $L^{1}(\mathbb{B}_n,\mathrm{d}\nu)$. Then the following statements hold true.
\begin{enumerate}
\item If for almost all $z\in\mathbb{B}_n$, $f(\mathrm{e}^{\mathrm{i}\gamma_1\theta}z_1,\ldots,\mathrm{e}^{\mathrm{i}\gamma_n\theta}z_n)=f(z)$ for almost all $\theta\in\mathbb{R}$, then whenever $l=(l_1,\ldots,l_n)\in\mathbb{Z}^n$ with $\gamma_1 l_1+\cdots+\gamma_n l_n\neq 0$, we have $\displaystyle\int\limits_{\mathbb{B}_n}f(z)z^{m+l}\bar{z}^{m}\mathrm{d}\nu(z)=0$ for all $m\in\mathbb{N}^n$ with $m+l\succeq 0$.
\item If the set $Z(f,l)=\displaystyle\{m\in\mathbb{N}^n: m+l\succeq 0\text{ and }\int\limits_{\mathbb{B}_n}f(z)z^{m+l}\bar{z}^{m}\mathrm{d}\nu=0\}$ does not have property (P) whenever $l=(l_1,\ldots,l_n)\in\mathbb{Z}^n$ with $\gamma_1 l_1+\cdots+\gamma_n l_n\neq 0$ then for almost all $z\in\mathbb{B}_n$, for almost all $\theta\in\mathbb{R}$, we have $f(\mathrm{e}^{\mathrm{i}\gamma_1\theta}z_1,\ldots,\mathrm{e}^{\mathrm{i}\gamma_n\theta}z_n)=f(z)$.
\end{enumerate}
\end{corollary}
\begin{proof}
Define $\displaystyle
g(z) = \dfrac{1}{2\pi}\int\limits_{0}^{2\pi}f(\mathrm{e}^{\mathrm{i}\gamma_1 t}z_1,\ldots,\mathrm{e}^{\mathrm{i}\gamma_n t}z_n)\mathrm{d}t,
$
for $z\in\mathbb{B}_n$ such that the integral on the right hand side is defined. Since $f\in L^{1}(\mathbb{B}_n,\mathrm{d}\nu)$, $g(z)$ is defined for almost all $z\in\mathbb{B}_n$ and for such $z$,  $g(\mathrm{e}^{\mathrm{i}\gamma_1\theta}z_1,\ldots,\mathrm{e}^{\mathrm{i}\gamma_n\theta}z_n)=g(z)$ for all $\theta\in\mathbb{R}$. Now for $l\in\mathbb{Z}^n$ and $m\in\mathbb{N}^n$ with $m+l\succeq 0$,
\begin{align*}
&\int\limits_{\mathbb{B}_n}g(z)z^{m+l}\bar{z}^m\mathrm{d}\nu(z)\\
&\quad\quad\quad =  \int\limits_{\mathbb{B}_n}\big\{\dfrac{1}{2\pi}\int\limits_{0}^{2\pi}f(\mathrm{e}^{\mathrm{i}\gamma_1 t}z_1,\ldots,\mathrm{e}^{\mathrm{i}\gamma_n t}z_n)\mathrm{d}t\big\}z^{m+l}\bar{z}^m\mathrm{d}\nu(z)\\
&\quad\quad\quad = \dfrac{1}{2\pi}\int\limits_{0}^{2\pi}\big\{\int\limits_{\mathbb{B}_n}f(\mathrm{e}^{\mathrm{i}\gamma_1 t}z_1,\ldots,\mathrm{e}^{\mathrm{i}\gamma_n t}z_n)z^{m+l}\bar{z}^m\mathrm{d}\nu\big\}\mathrm{d}t\\
&\quad\quad\quad = \dfrac{1}{2\pi}\int\limits_{0}^{2\pi}\big\{\int\limits_{\mathbb{B}_n}f(z_1,\ldots,z_n)z^{m+l}\bar{z}^m\mathrm{d}\nu\big\}\mathrm{e}^{-\mathrm{i}(\gamma_1 l_1+\cdots\gamma_n l_n)t}\mathrm{d}t\\
&\quad\quad\text{(by the invariance of the measure $\nu$ under the action of the $n-$torus)}\\
&\quad\quad\quad = \big(\dfrac{1}{2\pi}\int\limits_{0}^{2\pi}\mathrm{e}^{-\mathrm{i}(\gamma_1 l_1+\cdots+\gamma_n l_n)t}\mathrm{d}t\big)\big(\int\limits_{\mathbb{B}_n}f(z_1,\ldots,z_n)z^{m+l}\bar{z}^m\mathrm{d}\nu\big)\\
&\quad\quad\quad = \begin{cases}
0 & \text{if } \gamma_1 l_1 + \cdots + \gamma_n l_n \neq 0\\
\displaystyle\int\limits_{\mathbb{B}_n}f(z_1,\ldots,z_n)z^{m+l}\bar{z}^m\mathrm{d}\nu & \text{if } \gamma_1 l_1+\cdots+\gamma_n l_n = 0.
\end{cases}
\end{align*}

If for almost all $z\in\mathbb{B}_n$, for almost all $\theta\in\mathbb{R}$, $f(\mathrm{e}^{\mathrm{i}\gamma_1\theta}z_1,\ldots,\mathrm{e}^{\mathrm{i}\gamma_n\theta}z_n)=f(z)$ then $f(z)=g(z)$ for almost all $z\in\mathbb{B}_n$. The above computations then show that $\displaystyle\int\limits_{\mathbb{B}_n}f(z)z^{m+l}\bar{z}^{m}\mathrm{d}\nu(z)=0$ for all $m\in\mathbb{N}^n$ with $m+l\succeq 0$, whenever $\gamma_1 l_1 + \cdots + \gamma_n l_n \neq 0$.

Now suppose $Z(f,l)$ does not have property (P) whenever $\gamma_1 l_1 + \cdots + \gamma_n l_n \neq 0$. Then from the above computations, for all $l\in\mathbb{Z}^n$, the set of all $m\in\mathbb{N}^n$ with $m+l\succeq 0$ and $\displaystyle\int\limits_{\mathbb{B}_n}(f(z)-g(z))z^{m+l}\bar{z}^m\mathrm{d}\nu(z)=0$ does not have property (P). Corollary \ref{cor-1} now shows that $f(z)=g(z)$ for almost all $z\in\mathbb{B}_n$. Hence for almost all $z\in\mathbb{B}_n$, we have $f(\mathrm{e}^{\mathrm{i}\gamma_1\theta}z_1,\ldots,\mathrm{e}^{\mathrm{i}\gamma_n\theta}z_n)=f(z)$ for almost all $\theta\in\mathbb{R}$.
\end{proof}

For $\zeta=(\zeta_1,\ldots,\zeta_n)\in\mathbb{C}^n$ we write $\Sigma\zeta$ for $\zeta_1+\cdots+\zeta_n$. If $m=(m_1,\ldots,m_n)\in\mathbb{N}^n$ is a multi-index then we use the more common notation $|m|$ for $m_1+\cdots+m_n$ instead of $\Sigma m$.

For any bounded measurable function $g$ on $\mathbb{B}_n$, any $m\in\mathbb{N}^n$, and $\alpha>-1$, define
\begin{equation*}
\omega_{\alpha}(g,m) = \langle T_g e_{m}, e_{m}\rangle_{\alpha} = \int\limits_{\mathbb{B}_n}g(z)e_{m}(z)\bar{e}_{m}(z)\mathrm{d}\nu_{\alpha}(z).
\end{equation*}
The following theorem characterizes all $l\in\mathbb{Z}^n$ such that the set $\{m\in\mathbb{N}^n: m+l\succeq 0\text{ and } \omega_{\alpha}(g,m+l)=\omega_{\alpha}(g,m)\}$ does not have property (P), when $g$ has a special form.

\begin{proposition}\label{prop-2}
Suppose $g(z)=|z_1|^{2s_1}\cdots|z_n|^{2s_n}h(|z|)$ for $z\in\mathbb{B}_n$, where $s_1,\ldots, s_n\geq 0$ and $h:[0,1)\rightarrow\mathbb{C}$ is a bounded measurable function. Assume that $g$ is not a constant function on $\mathbb{B}_n$. Then for $l=(l_1,\ldots,l_n)\in\mathbb{Z}^n$ with $\Sigma l=0$ and $s_1l_1=\cdots=s_nl_n=0$, we have $\omega_{\alpha}(g,m+l)=\omega_{\alpha}(g,m)$ for all $m\in\mathbb{N}^n$ with $m+l\succeq 0$. Conversely, if $l=(l_1,\ldots,l_n)\in\mathbb{Z}^n$ such that the set $\{m\in\mathbb{N}^n: m+l\succeq 0 \text{ and } \omega_{\alpha}(g,m+l)=\omega_{\alpha}(g,m)\}$ does not have property (P) then $\Sigma l=0$ and $s_1 l_1=\cdots=s_nl_n=0$.
\end{proposition}
\begin{proof}
We first notice that when $n=1$, if $g(z_1)=|z_1|^{s_1}h(|z_1|)$ then we may rewrite $g(z_1)=\tilde{h}(|z_1|)$ where $\tilde{h}(t)=t^{s_1}h(t)$ for $0\leq t<1$. By this reason we always assume that $s_1=0$ if $n=1$.

We next recall the following formula. For any $\lambda=(\lambda_1,\ldots,\lambda_n)\in\mathbb{C}^n$ with $\Re(\lambda_1),\ldots,\Re(\lambda_n)>-1$, we have
\begin{equation}\label{eqn-3}
\int\limits_{\mathbb{S}_n}|\zeta_1|^{2\lambda_1}\cdots|\zeta_n|^{2\lambda_n}\mathrm{d}\sigma(\zeta) = \dfrac{\Gamma(n)\Gamma(\lambda_1+1)\cdots\Gamma(\lambda_n+1)}{\Gamma(n+\Sigma\lambda)}.
\end{equation}
The case $\lambda\in\mathbb{N}^n$ is in Lemma 1.11 in \cite{Zhu2005} and in fact the same argument works also for the general case.

Now for $\lambda=(\lambda_1,\ldots,\lambda_n)\in\mathbb{C}^n$ with $\Re(\lambda_1),\ldots,\Re(\lambda_n)>-1$, we have
\begin{align}
& \int\limits_{\mathbb{B}_n}|z_1|^{2\lambda_1}\cdots|z_n|^{2\lambda_n}h(|z|)\mathrm{d}\nu_{\alpha}(z)\notag\\
& =2nc_{\alpha}\int\limits_{0}^{1}r^{2n-1}\Big\{\int\limits_{\mathbb{S}_n}|r\zeta_1|^{2\lambda_1}\cdots|r\zeta_n|^{2\lambda_n}\mathrm{d}\sigma(\zeta)\Big\}h(r)(1-r^2)^{\alpha}\mathrm{d}r\notag\\
& = 2nc_{\alpha}\dfrac{\Gamma(n)\Gamma(\lambda_1+1)\cdots\Gamma(\lambda_n+1)}{\Gamma(n+\Sigma\lambda)}\int\limits_{0}^{1}r^{2n+2\Sigma\lambda-1}h(r)(1-r^2)^{\alpha}\mathrm{d}r\notag\\
& = c_{\alpha}\dfrac{\Gamma(n+1)\Gamma(\lambda_1+1)\cdots\Gamma(\lambda_n+1)}{\Gamma(n+\Sigma\lambda)}\int\limits_{0}^{1}r^{n+\Sigma\lambda-1}h(r^{1/2})(1-r)^{\alpha}\mathrm{d}r\notag\\
& = \dfrac{\Gamma(n+\alpha+1)}{\Gamma(\alpha+1)}\dfrac{\Gamma(\lambda_1+1)\cdots\Gamma(\lambda_n+1)}{\Gamma(n+\Sigma\lambda)}\int\limits_{0}^{1}r^{n+\Sigma\lambda-1}h(r^{1/2})(1-r)^{\alpha}\mathrm{d}r,\label{eqn-500}
\end{align}
since $c_{\alpha}=\dfrac{\Gamma(n+\alpha+1)}{\Gamma(n+1)\Gamma(\alpha+1)}$.

Put $s=(s_1,\ldots,s_n)$. For any $w\in\mathbb{C}$ with $\Re(w)\geq 1$, define
\begin{equation}\label{eqn-501}
H(w) = \dfrac{\Gamma(w+\alpha+1)}{\Gamma(w+\Sigma s)}\int\limits_{0}^{1}r^{w+\Sigma s-1}h(r^{1/2})(1-r)^{\alpha}\mathrm{d}r.
\end{equation}
Arguing as in the proof of Lemma \ref{lemma-1} we see that $H$ is analytic on the half plane $\Re(w)>1$ and is continuous on $\Re(w)\geq 1$. By the asymptotic behavior of the Gamma function at infinity and the boundedness of $h$, there is a polynomial $p$ such that $|H(w)|\leq p(|w|)$ for all $\Re(w)\geq 1$.

For $\zeta=(\zeta_1,\ldots,\zeta_n)\in\mathbb{C}^n$ with $\Re(\zeta_j)\geq 0$, define
\begin{align*}
F(\zeta) & = \dfrac{1}{\Gamma(\alpha+1)}\prod_{j=1}^{n}\dfrac{\Gamma(\zeta_j+s_j+1)}{\Gamma(\zeta_j+1)}H(n+\Sigma\zeta).
\end{align*}
Then $F$ is analytic in the interior of its defining domain and for any $m=(m_1,\ldots,m_n)$ in $\mathbb{N}^n$, we have
\begin{align*}
& \omega_{\alpha}(g,m)\\
& = \langle T_g e_m, e_m\rangle_{\alpha}\\
& = \dfrac{\Gamma(n+|m|+\alpha+1)}{\Gamma(n+\alpha+1)\prod_{j=1}^{n}\Gamma(m_j+1)}\int\limits_{\mathbb{B}_n}g(z)z^{m}\bar{z}^{m}\mathrm{d}\nu_{\alpha}(z)\\
& = \dfrac{\Gamma(n+|m|+\alpha+1)}{\Gamma(n+\alpha+1)\prod_{j=1}^{n}\Gamma(m_j+1)}\int\limits_{\mathbb{B}_n}|z_1|^{2(s_1+m_1)}\cdots|z_n|^{2(s_n+m_n)}h(|z|)\mathrm{d}\nu_{\alpha}(z)\\
& = \dfrac{1}{\Gamma(\alpha+1)}\prod_{j=1}^{n}\dfrac{\Gamma(m_j+s_j+1)}{\Gamma(m_j+1)}H(n+|m|)\quad\quad\text{(by \eqref{eqn-500} and \eqref{eqn-501})}\\
& = F(m).
\end{align*}

Suppose $l=(l_1,\ldots,l_n)\in\mathbb{Z}^n$ such that $s_1l_1=\cdots=s_nl_n=0$ and $\Sigma l=0$. Then for all $m\in\mathbb{N}^n$ with $m+l\succeq 0$ we have $\dfrac{\Gamma(m_j+s_j+1)}{\Gamma(m_j+1)}=\dfrac{\Gamma(m_j+l_j+s_j+1)}{\Gamma(m_j+l_j+1)}$ for $1\leq j\leq n$, and $|m+l|=|m|+\Sigma l=|m|$. Hence we have $\omega_{\alpha}(g,m+l)=F(m+l)=F(m)=\omega_{\alpha}(g,m)$ for all such $m$.

Conversely, suppose $l=(l_1,\ldots,l_n)\in\mathbb{Z}^n$ such that the set $\{m\in\mathbb{N}^n: m+l\succeq 0 \text{ and } \omega_{\alpha}(g,m+l)=\omega_{\alpha}(g,m)\}$ does not have property (P). Then the set $\{m\in\mathbb{N}^n: m+l\succeq 0 \text{ and } F(m+l)=F(m)\}$ does not have property (P). By Proposition \ref{prop-1}, we see that $F(\zeta+l)=F(\zeta)$ for all $\zeta\in\mathbb{C}^n$ with $\Re(\zeta_j)>\max\{0,-l_j\}, j=1,\ldots, n$. This implies that for such $\zeta$,
\begin{equation}\label{eqn-301}
\prod_{j=1}^{n}{\dfrac{\Gamma(\zeta_j+s_j+l_j+1)}{\Gamma(\zeta_j+l_j+1)}}H(n+\Sigma\zeta+\Sigma l)
= \prod_{j=1}^{n}{\dfrac{\Gamma(\zeta_j+s_j+1)}{\Gamma(\zeta_j+1)}}H(n+\Sigma\zeta).
\end{equation}

We now show that $s_1 l_1=0$. If $s_1=0$, there is nothing to show. So suppose $s_1>0$ (then we must have $n\geq 2$). We will show that $l_1=0$. Assume for contradiction that $l_1\neq 0$. We consider here only the case $l_1>0$. The case $l_1<0$ can be handled is a similar fashion. Since $l_1\geq 1$, equation \eqref{eqn-301} gives
\begin{align}\label{eqn-302}
& (\zeta_1+s_1+1)\cdots(\zeta_1+s_1+l_1)\prod_{j=2}^{n}\dfrac{\Gamma(\zeta_j+s_j+l_j+1)}{\Gamma(\zeta_j+l_j+1)} H(n+\Sigma\zeta+\Sigma l)\notag\\
&= (\zeta_1+1)\cdots (\zeta_1+l_1)\prod_{j=2}^{n}\dfrac{\Gamma(\zeta_j+s_j+1)}{\Gamma(\zeta_j+1)}H(n+\Sigma\zeta).
\end{align}
Now choose $\tilde{m}=(m_2,\ldots, m_n)$ such that $m_j\geq\max\{1,1-l_j\}$ for all $j=2,\ldots,m$ and such that $\displaystyle\int\limits_{0}^{1}r^{n+|\tilde{m}|+\Sigma l+\Sigma s-2}h(r^{1/2})(1-r)^{\alpha}\mathrm{d}r\neq 0$. It is possible to do so since $h$ is not identically zero on $[0,1)$. Then with $\zeta=(\zeta_1,\tilde{m})$, \eqref{eqn-302} gives
\begin{align*}
& (\zeta_1+s_1+1)\cdots(\zeta_1+s_1+l_1)\prod_{j=2}^{n}\dfrac{\Gamma(m_j+s_j+l_j+1)}{\Gamma(m_j+l_j+1)} H(n+\zeta_1+|\tilde{m}|+\Sigma l)\\
&= (\zeta_1+1)\cdots (\zeta_1+l_1)\prod_{j=2}^{n}\dfrac{\Gamma(m_j+s_j+1)}{\Gamma(m_j+1)} H(n+\zeta_1+|\tilde{m}|),
\end{align*}
for all $\zeta_1\in\mathbb{C}^n$ with $\Re(\zeta_1)>0$. Since each side of the above identity is in fact an analytic function of $\zeta_1$ on $\Re(\zeta_1)>\max\{1-n-|\tilde{m}|-\Sigma l,1-n-|\tilde{m}|\}$ (which contains $-1$), the identity still holds true for $\zeta_1=-1$. Therefore, we get
\begin{equation*}
s_1(s_1+1)\cdots (s_1+l_1-1)\prod_{j=2}^{n}\dfrac{\Gamma(m_j+s_j+l_j+1)}{\Gamma(m_j+l_j+1)} H(n+|\tilde{m}|+\Sigma l-1) = 0,
\end{equation*}
which is a contradiction because the left hand side is nonzero by the choice of $m_2,\ldots,m_n$. Thus we have $l_1=0$. Similarly, we have $s_j l_j=0$ for any $j=2,\ldots,n$. Now equation \eqref{eqn-301} becomes $
H(n+\Sigma\zeta+\Sigma l)= H(n+\Sigma\zeta)$ for all $\zeta\in\mathbb{C}^n$ with $\Re(\zeta_j)>\max\{0,-l_j\}, j=1,\ldots,n$. Let $w=n+\Sigma\zeta$. Then $\Re(w)>\max\{n,n-\Sigma l\}$ and $H(w+\Sigma l)=H(w)$.

We now show that $\Sigma l=0$. Assume for contradiction that $\Sigma l\neq 0$. By changing $w$ to $w+\Sigma l$ if necessary, we may assume that $\Sigma l>0$. Since $H$ is periodic and analytic in $\Re(w)>n$, it extends to a periodic entire function on $\mathbb{C}$ which we still denote by $H$. Now for a complex number $w$ with $\Re(w)\leq n$, choose a number $k\in\mathbb{N}$ such that $\Re(w)+k(\Sigma l)\leq n< \Re(w)+(k+1)(\Sigma l)$. Then we have
\begin{align*}
|H(w)| & = |H(w+(k+1)(\Sigma l))|\leq p(|w+k(\Sigma l)+\Sigma l|)\leq q(|w+k(\Sigma l)|),
\end{align*}
for some polynomial $q$ with nonnegative coefficients. Now since $k(\Sigma l)<n-\Re(w)\leq n+|w|$, we have $q(|w+k(\Sigma l)|)\leq q(|w|+k(\Sigma l))\leq q(n+2|w|)$. Hence $|H(w)|\leq q(n+2|w|)$ when $\Re(w)\leq n$. Now when $\Re(w)>n$, we have $|H(w)|\leq p(|w|)$. In general, for any $w\in\mathbb{C}$, $|H(w)|\leq q(n+2|w|)+p(|w|)$. Since $H$ is entire, it must be a polynomial. But $H$ is also periodic, so it must be a constant function. We conclude that there is a constant $c$ so that
\begin{equation*}
\dfrac{\Gamma(w+\alpha+1)}{\Gamma(w+\Sigma s)}\int\limits_{0}^{1}r^{w+\Sigma s-1}h(r^{1/2})(1-r)^{\alpha}\mathrm{d}r=H(w) = c
\end{equation*}
for all $w\in\mathbb{C}$ with $\Re(w)\geq 1$.

Suppose $\Sigma s=0$ (hence $s_1=\cdots=s_n=0$). Then from the identity
\begin{equation*}
1 = \dfrac{\Gamma(w+\alpha+1)}{\Gamma(w)\Gamma(\alpha+1)}\int\limits_{0}^{1}r^{w-1}(1-r)^{\alpha}\mathrm{d}r,
\end{equation*}
we conclude that $\displaystyle\int\limits_{0}^{1}\Big(h(r^{1/2})-\dfrac{c}{\Gamma(\alpha+1)}\Big)r^{w-1}(1-r)^{\alpha}\mathrm{d}r=0$ for all $w\in\mathbb{C}$ with $\Re(w)\geq 1$. Thus, $h(t)=\dfrac{c}{\Gamma(\alpha+1)}$ for almost all $t\in [0,1)$. This shows that $g(z)=h(|z|)$ is a constant function on $\mathbb{B}_n$, which is a contradiction.

Now suppose $\Sigma s>0$. Then for $w=u\in\mathbb{R}$ with $u\geq 1$, we have
\begin{align*}
|c| & = \Big|\dfrac{\Gamma(u+\alpha+1)}{\Gamma(u+\Sigma s)}\int\limits_{0}^{1}r^{u+\Sigma s-1}h(r^{1/2})(1-r)^{\alpha}\mathrm{d}r\Big|\\
& \leq \|h\|_{\infty}\dfrac{\Gamma(u+\alpha+1)}{\Gamma(u+\Sigma s)}\int\limits_{0}^{1}r^{u+\Sigma s-1}(1-r)^{\alpha}\mathrm{d}r\\
& = \|h\|_{\infty}\dfrac{\Gamma(u+\alpha+1)}{\Gamma(u+\Sigma s)}\dfrac{\Gamma(\alpha+1)\Gamma(u+\Sigma s)}{\Gamma(u+\Sigma s+\alpha+1)}\\
& = \|h\|_{\infty}\dfrac{\Gamma(\alpha+1)\Gamma(u+\alpha+1)}{\Gamma(u+\Sigma s+\alpha+1)}\\
& \approx \|h\|_{\infty}\Gamma(\alpha+1)u^{-\Sigma s},
\end{align*}
by Stirling's formula for the Gamma function. Let $u\rightarrow\infty$, we get $c=0$. So we have $\displaystyle\int\limits_{0}^{1}r^{w+\Sigma s-1}h(r^{1/2})(1-r)^{\alpha}\mathrm{d}r = 0$ for all $w\in\mathbb{C}$ with $\Re(w)\geq 1$. This implies that $h(r)=0$ for almost all $r\in [0,1)$. Hence $g(z)=0$ for almost all $z\in\mathbb{B}_n$, which is again a contradiction. Thus we have $\Sigma l=0$.
\end{proof}

We are now ready for the proof of Theorem \ref{theorem-2}.

\begin{proof}[Proof of Theorem \ref{theorem-2}]
Since $g(z_1,\ldots,z_n)=g(|z_1|,\ldots,|z_n|)$ for almost all $z\in\mathbb{B}_n$, Theorem 3.1 in \cite{Le-5} shows that the Toeplitz operator $T_g$ is diagonal with respect to the standard orthonormal basis. The eigenvalues of $T_g$ are given by $\omega_{\alpha}(g,m)=\langle T_g e_m, e_m\rangle_{\alpha}$ for $m\in\mathbb{N}^n$. Note that $\omega_{\alpha}(\bar{g},m)=\bar{\omega}_{\alpha}(g,m)$ for all such $m$.

Now $T_fT_g=T_gT_f$ on $\mathcal{P}$ if and only if for all $l\in\mathbb{Z}^n$ and $m\in\mathbb{N}^n$ with $m+l\succeq 0$,
\begin{align}\label{eqn-304}
0 & = \langle (T_fT_g-T_gT_f)e_{m+l}, e_{m}\rangle_{\alpha}\notag\\
& = \langle \omega_{\alpha}(g,m+l) T_f e_{m+l}, e_{m}\rangle_{\alpha} - \langle T_f e_{m+l}, \omega_{\alpha}(\bar{g},m)e_{m}\rangle_{\alpha}\\
& = (\omega_{\alpha}(g,m+l)-\omega_{\alpha}(g,m))\langle T_f e_{m+l}, e_{m}\rangle_{\alpha}.\notag
\end{align}

Suppose $f(\mathrm{e}^{\mathrm{i}\theta}z)=f(z)$ for almost all $\theta\in\mathbb{R}$, almost all $z\in\mathbb{B}_n$, and for $1\leq j\leq n$ with $s_j\neq 0$, $f(z_1,\ldots,z_{j-1},|z_j|,z_{j+1},\ldots,z_n)=f(z)$ for almost all $z\in\mathbb{B}_n$. Let $l=(l_1,\ldots,l_n)$ be in $\mathbb{Z}^n$. If $\Sigma l\neq 0$ or for some $1\leq j\leq n$, $s_jl_j\neq 0$ (hence $s_j\neq 0$ and $l_j\neq 0$) then conclusion (1) of Corollary \ref{cor-2} shows that $\langle T_f e_{m+l}, e_{m}\rangle_{\alpha}=0$ for all $m\in\mathbb{N}^n$ with $m+l\succeq 0$. If $\Sigma l=0$ and $s_1l_1=\cdots=s_n l_n=0$, Proposition \ref{prop-2} shows that $\omega_{\alpha}(g,m+l)=\omega_{\alpha}(g,m)$ for all $m\in\mathbb{N}^n$ with $m+l\succeq 0$. Thus \eqref{eqn-304} holds for all $l\in\mathbb{Z}^n$ and $m\in\mathbb{N}^n$ with $m+l\succeq 0$. Therefore $T_fT_g=T_gT_f$ on $\mathcal{P}$.

Now suppose $T_fT_g=T_gT_f$ on $\mathcal{P}$. Let $l=(l_1,\ldots,l_n)$ be in $\mathbb{Z}^n$ such that $\Sigma l\neq 0$ or $s_j l_j\neq 0$ for some $1\leq j\leq n$. Then Proposition \ref{prop-2} shows that the set $\{m\in\mathbb{N}^n: m+l\succeq 0\text{ and } \omega_{\alpha}(g,m+l)=\omega_{\alpha}(g,m)\}$ has property (P). Since \eqref{eqn-304} holds for all $m\in\mathbb{N}^n$ with $m+l\succeq 0$, we conclude that the set $\displaystyle\{m\in\mathbb{N}^n: m+l\succeq 0 \text{ and }\int\limits_{\mathbb{B}_n}f(z)z^{m+l}\bar{z}^{m}(1-|z|^2)^{\alpha}\mathrm{d}\nu(z)=0\}$ does not have property (P). This is true whenever $l\in\mathbb{Z}^n$ such that $\Sigma l\neq 0$ or $s_j l_j\neq 0$ for some $1\leq j\leq n$. Conclusion (2) of Corollary \ref{cor-2} now implies that for almost all $\theta\in\mathbb{R}$ and almost all $z\in\mathbb{B}_n$, we have $f(z)=f(\mathrm{e}^{\mathrm{i}\theta}z)=f(z_1,\ldots,z_{j-1},\mathrm{e}^{\mathrm{i}\theta}z_{j},z_{j+1},\ldots,z_n)$ for any $1\leq j\leq n$ with $s_j\neq 0$. This shows that $f(\mathrm{e}^{\mathrm{i}\theta}z)=f(z)$ for almost all $\theta\in\mathbb{R}$, almost all $z\in\mathbb{B}_n$, and for $1\leq j\leq n$ with $s_j\neq 0$, $f(z_1,\ldots,z_{j-1},|z_j|,z_{j+1},\ldots,z_n)=f(z)$ for almost all $z\in\mathbb{B}_n$.
\end{proof}

\begin{remark}
If $s_1=\cdots=s_n=0$ so that $g(z)=h(|z|)$ for a non-constant bounded measurable function $h$ on $[0,1)$, then Theorem \ref{theorem-2} shows that for $f\in L^2_{\alpha}$, $T_f$ commutes with $T_g$ if and only if $f(\mathrm{e}^{\mathrm{i}\theta}z)=f(z)$ for almost all $z\in\mathbb{B}_n$, almost all $\theta\in\mathbb{R}$. In the one dimensional case, those functions are exactly radial functions. So we recover {\v{C}}u{\v{c}}kovi{\'c} and Rao's result.
\end{remark}

\begin{remark}
If $g(z)=|z_1|\cdots |z_{n-1}|h(|z|)$ for some bounded measurable function $h$ on $[0,1)$ then Theorem \ref{theorem-2} shows that for $f\in L^{2}_{\alpha}$, $T_f$ commutes with $T_g$ if and only if $f(z)=f(\mathrm{e}^{\mathrm{i}\theta}z)=f(|z_1|,\ldots,|z_{n-1}|,z_n)$ for almost all $\theta\in\mathbb{R}$ and almost all $z\in\mathbb{B}_{n}$ . This is equivalent to the condition that $f(z)=f(|z_1|,\ldots,|z_n|)$ for almost all $z\in\mathbb{B}_n$.
\end{remark}


\begin{thebibliography}{10}

\bibitem{Axler1991}
Sheldon Axler and {\v{Z}}eljko {\v{C}}u{\v{c}}kovi{\'c}, \emph{Commuting
  {T}oeplitz operators with harmonic symbols}, Integral Equations Operator
  Theory \textbf{14} (1991), no.~1, 1--12. \MR{MR1079815 (92f:47018)}

\bibitem{Axler2000}
Sheldon Axler, {\v{Z}}eljko {\v{C}}u{\v{c}}kovi{\'c}, and N.~V. Rao,
  \emph{Commutants of analytic {T}oeplitz operators on the {B}ergman space},
  Proc. Amer. Math. Soc. \textbf{128} (2000), no.~7, 1951--1953. \MR{MR1694299
  (2000m:47035)}

\bibitem{Brown1963}
Arlen Brown and Paul~R. Halmos, \emph{Algebraic properties of {T}oeplitz
  operators}, J. Reine Angew. Math. \textbf{213} (1963/1964), 89--102.
  \MR{MR0160136 (28 \#3350)}

\bibitem{Coburn1973}
Lewis~A. Coburn, \emph{Singular integral operators and {T}oeplitz operators on
  odd spheres}, Indiana Univ. Math. J. \textbf{23} (1973/74), 433--439.
  \MR{MR0322595 (48 \#957)}

\bibitem{Cuckovic1994}
{\v{Z}}eljko {\v{C}}u{\v{c}}kovi{\'c}, \emph{Commutants of {T}oeplitz operators
  on the {B}ergman space}, Pacific J. Math. \textbf{162} (1994), no.~2,
  277--285. \MR{MR1251902 (94j:47041)}

\bibitem{Cuckovic1998}
{\u{Z}}eljko {\u{C}}u{\u{c}}kovi{\'c} and N.~V. Rao, \emph{Mellin transform,
  monomial symbols, and commuting {T}oeplitz operators}, J. Funct. Anal.
  \textbf{154} (1998), no.~1, 195--214. \MR{MR1616532 (99f:47033)}

\bibitem{Le-5}
Trieu Le, \emph{Diagonal {T}oeplitz operators on weighted {B}ergman spaces},
  preprint.

\bibitem{Rudin1987}
Walter Rudin, \emph{Real and complex analysis}, third ed., McGraw-Hill Book
  Co., New York, 1987. \MR{MR924157 (88k:00002)}

\bibitem{Zhu2005}
Kehe Zhu, \emph{Spaces of holomorphic functions in the unit ball}, Graduate
  Texts in Mathematics, vol. 226, Springer-Verlag, New York, 2005.
  \MR{MR2115155 (2006d:46035)}

\bibitem{Zhu2007}
\bysame, \emph{Operator theory in function spaces}, second ed., Mathematical
  Surveys and Monographs, vol. 138, American Mathematical Society, Providence,
  RI, 2007. \MR{MR2311536}

\end{thebibliography}
\end{document}